\documentclass[oneside,english]{amsart}
\usepackage[T1]{fontenc}
\usepackage[latin9]{inputenc}
\usepackage{verbatim}
\usepackage{float}
\usepackage{amsbsy}
\usepackage{amstext}
\usepackage{amsthm}
\usepackage{amssymb}
\usepackage{graphicx}
\PassOptionsToPackage{normalem}{ulem}
\usepackage{ulem}
 \pdfoutput=1

\makeatletter

\floatstyle{ruled}
\newfloat{algorithm}{tbp}{loa}
\providecommand{\algorithmname}{Algorithm}
\floatname{algorithm}{\protect\algorithmname}

\numberwithin{equation}{section}
\numberwithin{figure}{section}
\theoremstyle{plain}
\newtheorem{thm}{\protect\theoremname}[section]
\theoremstyle{remark}
\newtheorem{rem}[thm]{\protect\remarkname}
\theoremstyle{plain}
\newtheorem{cor}[thm]{\protect\corollaryname}
\theoremstyle{plain}
\newtheorem{lem}[thm]{\protect\lemmaname}

\usepackage{a4wide}

\usepackage{bbm}
\usepackage{enumerate}
\usepackage{tikz}

\allowdisplaybreaks

\newcommand{\R}{\mathbb{R}}
\newcommand{\p}{\mathbb{P}}
\newcommand{\E}{\mathbb{E}}

\newcommand{\N}{\mathbb{N}}

\newcommand{\1}{\mathbbm{1}}

\newcommand{\convN}{\overset{\text{a.s.}}{\underset{N\rightarrow+\infty}{\longrightarrow}}}

\colorlet{mycol}{black}
\newcommand\mybox[3][]{\tikz[anchor=base,baseline]\node[inner sep=2pt,draw=#2,#1]{$\displaystyle#3\mathstrut$};}

\newtheorem{hypothesis}{Hypothesis}

\date{\today}

\AtBeginDocument{
  
}

\makeatother

\usepackage{babel}
\usepackage{listings}
\providecommand{\corollaryname}{Corollary}
\providecommand{\lemmaname}{Lemma}
\providecommand{\remarkname}{Remark}
\providecommand{\theoremname}{Theorem}

\begin{document}
\title[Counterexample to a transition probability formula]{Counterexample to a transition probability formula for the ancestral
process}
\author{Sylvain Rubenthaler}
\begin{abstract}
We consider weighted particle systems in which new generations are
re-sampled from current particles with probabilities proportional
to their weights. This covers a broad class of sequential Monte Carlo
methods, widely used in applied statistics. We consider the genealogical
tree embedded into such particle system. When the time is reversed,
the particle system induces a partition valued family of processes
(partitions on the leaves of the genealogical tree). Our aim here
is to give a counterexample to a well known formula describing the
transition probabilities of this process. 
\end{abstract}

\maketitle

\section{Introduction}

\subsection{Description of the problem}

We consider interacting particle systems (IPS). These are a broad
class of stochastic models for phenomena in disciplines including
physics, engineering, biology and finance. Important examples are
particle filters, particle methods and sequential Monte Carlo (SMC),
which feature prominently in numerical approximation schemes for nonlinear
filtering, as well as mean field approximation of Feynman-Kac flows.
For additional background, we direct readers to \cite{del-moral-2004},
\cite{crisan-rozovski-2011}.

Central to these methods are discrete time, evolving weighted particle
systems. Correlations between particles arise out of resampling: a
stochastic selection mechanism in which particles with high weight
are typically replicated while those with low weight vanish, giving
rise to an embedded genealogy. Many papers study this genealogy. The
standard way is to reverse the time: we start with $n$ particles,
they have a certain number of ancestors (less than $n$) and so on.
One can carry on this way up to the first common ancestor of all the
particles. At each step, the genealogical structure can be represented
with a partition on our $n$ particles at the start: namely, two particles
are in the same block of they have the same ancestor at the given
step. Some authors rely on a transition probability formula to describe
the ancestral process on partitions (Equation (2) in \cite{mohle-1998},
Equation (1) in \cite{koskela-jenkins-johansen-spano-2020}). The
contribution of this paper is to give a counterexample to this formula
in the SMC case. 

\subsection{Description of the mathematical model. }

We consider a sub-case of the model described in \cite{koskela-jenkins-johansen-spano-2020}.
We have a measurable space $(E,\mathcal{E})$, a probability measure
$\mu$ on $(E,\mathcal{E})$, a Markov transition $K:E\times\mathcal{E}\rightarrow(0,\infty)$,
and a non-negative potential $g:E\rightarrow(0,\infty)$. These correspond
to the state space, initial proposal distribution, transition kernel,
and importance weight function of our IPS, respectively. 

Let $\zeta_{t}^{(N)}:=\{(w_{t}^{(i)},X_{t}^{(i)})\}_{i=1}^{N}$be
a weighted $N$-particle system at time $t$ in $\N$, where each
$X_{t}^{(i)}\in E$, and the weights $w_{t}^{(i)}$are non-negative
and satisfy $\sum_{i=1}^{N}w_{t}^{(i)}=1$. Let $S$ be a resampling
operator which acts on $\zeta_{t}^{(N)}$by assigning to each particle
a random number of offsprings. The total number of offsprings is fixed
at $N$, and the mean number of offsprings of particle $i\in\{1,\dots,N\}$
is $Nw_{t}^{(i)}$. All offsprings are assigned an equal weight $1/N$.
More concretely
\[
S\zeta_{t}^{(N)}=\{(N^{-1},X_{t}^{(a_{t}^{(i)})})\}_{i=1}^{N}\,,
\]
where $a_{t}^{(i)}=j$ if $j$ in $\zeta_{t}^{(N)}$ is the parent
of $i$ in $S\zeta_{t}^{(N)}$. Particles with low weight are randomly
removed by having no offspring, while particles with high weights
tend to have many offsprings.

The step from time $t$ to time $t+1$ is completed by propagating
each particle in the $N$-uple $S\zeta_{t}^{(N)}$independently through
the transition kernel $K$ to obtain particle locations $X_{t+1}^{(i)}\sim K(X_{t}^{(a_{t}^{(i)})},.)$.
Finally, each particle $i\in\{1,\dots,N\}$ is assigned a weight proportional
to the potential $g$ evaluated at the location of the particle, so
that the full update is
\[
\zeta_{t}^{(N)}\overset{S}{\mapsto}\{N^{-1},X_{t}^{(a_{t}^{(i)})}\}_{i=1}^{N}\overset{K}{\mapsto}\left\{ \frac{g(X_{t+1}^{(i)})}{\sum_{j=1}^{N}g(X_{t+1}^{(j)})},X_{t+1}^{(i)}\right\} _{i=1}^{N}\,.
\]
A specification is given in Algorithm \ref{alg:Simulation-of-an}.
\begin{algorithm}[h]
\begin{raggedright}
\texttt{1}:\texttt{~for} $i\in\{1,\dots,N\}$ \texttt{do}~\\
\texttt{2:~~~Sample $X_{0}^{(i)}\sim\mu$.}~\\
\texttt{3:~for $i\in\{1,\dots,N\}$ do}~\\
\texttt{4:~~~Set $w_{0}^{(i)}\leftarrow\frac{g(X_{0}^{(i)})}{g(X_{0}^{(1)})+\dots+g(X_{0}^{(N)})}$.}~\\
\texttt{5:~for $t\in\{0,\dots,T-1\}$ do}~\\
\texttt{6:~~~Sample $(a_{t}^{(1)},\dots,a_{t}^{(N)})\sim\text{Categorical}(w_{t}^{(1)},\dots,w_{t}^{(N)})$.}~\\
\texttt{7:~~~for $i\in\{1,\dots,N\}$ do}~\\
\texttt{8:~~~~~Sample $X_{t+1}^{(i)}\sim K(X_{t}^{(a_{t}^{(i)})},.)$.}~\\
\texttt{9:~~~for $i\in\{1,\dots,N\}$ do}~\\
\texttt{10:~~~~Set $w_{t+1}^{(i)}\leftarrow\frac{g(X_{t+1}^{(i)})}{g(X_{t+1}^{(1)})+\dots+g(X_{t+1}^{(N)})}$}
\par\end{raggedright}
\caption{Simulation of an IPS\label{alg:Simulation-of-an}}

\end{algorithm}
There is a genealogy embedded in Algorithm \ref{alg:Simulation-of-an}.
Consider $\zeta_{t}^{(N)}$ at a fixed time $t$. Tracing the ancestor
indices $(a_{t}^{(1)},\dots,a_{t}^{(N)})$ backward in time results
in a coalescing forest of lineages. The forest forms a tree once the
most recent common ancestor (MRCA) of all particles is reached, provided
that happens before reaching the initial time $0$. In \cite{koskela-jenkins-johansen-spano-2020},
the authors show that, under certain conditions and an appropriate
time-rescaling functionals of this genealogy depending upon finite
numbers of leaves ($n$) converge to corresponding functionals of
the Kingman $n$-coalescent \cite{kingman-1982} as $N\rightarrow\infty$
in the sense of finite dimensional distributions. In particular, they
show that the expected number of time steps from the leaves to the
MRCA scales linearly in $N$ for any finite number of leaves. Other
papers estimate the number of time steps from the leaves to the MRCA
by simpler methods (for example: \cite{del-moral-miclo-patras-rubenthaler-2010,jacob-murray-rubenthaler-2015}).

The rest of this paper is structured as follows. In Section \ref{sec:The-genealogical-process},
we present the genealogical process and its transition formula. In
Section \ref{sec:Counterexample}, we present our counterexample.
The main results are Theorem \ref{th:main} and Corollary \ref{cor:For-well-chosen}.
We conclude this section by summarizing notations. 

Let $(x)_{b}:=x(x-1)\dots(x-b+1)$ be the falling factorial. We adopt
the convention $\sum_{\emptyset}=0$, $\prod_{\emptyset}=1$. For
an integer $n\in\N$, we define $[n]:=\{1,2,\dots,n\}$ with $[0]:=\emptyset$.
For a partition $\xi$, $|\xi|$ denotes the number of blocks in $\xi$.%
\begin{comment}
and for a finite set $A$, we let $\Pi_{n}(A)$ denote the set of
partitions of $A$ into at most $n$ nonempty blocks, with $\Pi_{n}(\emptyset)=(\emptyset,\dots,\emptyset)$.
For a partition $\xi$, $|\xi|$ denotes the number of blocks in $\xi$,
and $\boldsymbol{x}:=(x_{1},\dots,x_{N})$, where the length of the
vector will be clear from context. For a vector $\boldsymbol{x}$,
$|\boldsymbol{x}|$ denotes the $L^{1}$-norm.
\end{comment}
{} 

\section{The genealogical process and its transition formula\label{sec:The-genealogical-process}}

It will be convenient to express our IPS in reverse time by denoting
the initial time in Algorithm \ref{alg:Simulation-of-an} by $T$
and the terminal time by $0$, and to describe the genealogy in terms
of a partition valued family of processes $(G_{t}^{(n,N)})_{t=0}^{T}$
indexed by $n\leq N$, where $n$ denotes the number of observed leaves
(time $0$ particles) in a system with $N$ particles. The process
$(G_{t}^{(n,N)})_{t=0}^{T}$ is defined in terms of the underlying
IPS via its initial condition $G_{0}^{(n,N)}=\{\{1\},\dots,\{n\}\},$
and its dynamics, which are driven by the requirement that $i\neq j\in[n]$
belong to the same block in $G_{t}^{(n,N)}$ if leaves $i$ and $j$
have a common ancestor at time $t$.

As said in \cite{koskela-jenkins-johansen-spano-2020}, the genealogical
process $(G_{t}^{(n,N)})_{t=0}^{T}$ evolves on a space which tracks
the ancestral relationships of the observed particles but not their
states. The process is a projection of the time reversal of the historical
process of \cite{del-moral-miclo-2001}, in which particle location
have been marginalized over. A consequence of this is that $(G_{t}^{(n,N)})_{t=0}^{T}$
is not Markovian in general.

Let $\nu_{t}^{(i)}$ denote the number of offsprings that particle
$i$ at time $t$ has at time $t-1$ (with $\boldsymbol{\nu}_{t}=(\nu_{t}^{(i)})_{1\leq i\leq N}$).
Let $\boldsymbol{a}_{t}=(a_{t}^{(i)})_{1\leq i\leq N}$ for all $t$.
We make the same assumption as in \cite{koskela-jenkins-johansen-spano-2020}
(``Standing assumption'', p. 5).

\begin{hypothesis} \label{hyp:parental-distribution}The conditional
distribution of parental indices given offsprings counts, $\boldsymbol{a}_{t}|\boldsymbol{\nu}_{t}$,
is uniform over all vectors which satisfy $\nu_{t}^{(i)}=\#\{j\in[N]\,:\,a_{t}^{(j)}=i\}$
for each $i\in[N]$.

\end{hypothesis}
\begin{rem}
The above assumption concerns the marginal distribution of parental
assignments without particle locations. A sufficient condition for
this assumption is exchangeability of the \texttt{Categorical} variables
in line 6 of Algorithm \ref{sec:The-genealogical-process} (this is
indeed the case). 

In \cite{mohle-1998}, there is the following additional assumption
(p. 439). 

\begin{hypothesis}\label{hyp:nu-independants}The offspring counts
($\boldsymbol{\nu}_{0}$, $\boldsymbol{\nu}_{1}$, ...) are independent.

\end{hypothesis}

We do not make this assumption since it is not likely to be true for
our IPS. 

Let $\xi$ and $\eta$ be partitions of $[n]$ with blocks ordered
by the least element in each block, and with $\eta$ obtained from
$\xi$ by merging some subsets of blocks. For $i\in[|\eta|],$let
$b_{i}$ be the number of blocks of $\xi$ that have been merged to
form block $i$ in $\eta$, so $b_{1}+\dots+b_{|\eta|}=|\xi|$, and
define
\begin{equation}
p_{\xi,\eta}(t):=\frac{1}{(N)_{|\xi|}}\sum_{\begin{array}{c}
i_{1}\neq\dots\neq i_{|\eta|}=1\\
\text{all distincts}
\end{array}}(\nu_{t}^{(i_{1})})_{b_{1}}\dots(\nu_{t}^{(i_{|\nu|})})_{b_{|\nu|}}\label{eq:transition-formula}
\end{equation}
(this is Equation (1) of \cite{koskela-jenkins-johansen-spano-2020}
and Equation (2) of \cite{mohle-1998}). \cite{koskela-jenkins-johansen-spano-2020,mohle-1998}
claim that the above is the transition probability from state $\xi$
at time $t-1$ to state $\eta$ at time $t$, given $\boldsymbol{\nu}_{t}$
(suppressed from the notation). The argument to support this is the
same in the two articles (\cite{koskela-jenkins-johansen-spano-2020,mohle-1998}):
we can think about proceeding from generation $t-1$ to $t$ as an
experiment in which $|\xi|$ balls (the children) are thrown into
$|\nu|$ boxes, such that at the end of the experiment the $i$-th
box contains exactly $\nu_{t}^{(i)}$ balls. 

\begin{comment}
There is a first underlying assumption in this reasoning. We are interested
in the black and white particle because it is the ancestor of a partition
block. If it chooses the same ancestor at time $t$ than another ancestor
of a partition block, then these two blocks merge. 
\end{comment}

There is an \uline{underlying assumption} to this reasoning in
\cite{koskela-jenkins-johansen-spano-2020,mohle-1998} and we are
going to expose it here. Suppose we have a partition $\xi$ at time
$t-1$. In Figure \ref{fig:Choice-of-the}, we have the genealogical
structure of the particle process. 
\begin{figure}[h]
\begin{centering}
\includegraphics{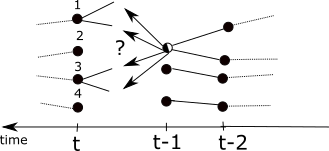}
\par\end{centering}
\caption{Choice of the ancestor\label{fig:Choice-of-the}}

\end{figure}
The partition $\xi$ is represented by three particles (for each of
them, its descendants at time $0$ form a block of $\xi$). We want
to make a computation conditioned to the knowledge of $\boldsymbol{\nu}_{t}$,
the knowledge of $\xi$%
\begin{comment}
and the knowledge of the \foreignlanguage{american}{representatives}
of each block of $\xi$ at time $t-1$ (i.e. the common ancestors
of the leaves in a block of $\xi$)
\end{comment}
. In our case, we know that two particles at time $t$ have two offsprings
and two particles at time $t$ have zero offsprings. This is represented
in our figure by ``dangling'' branches. Now, suppose we begin by
choosing an parent for the particle in black and white (which is the
common ancestor of leaves forming a block of $\xi$). We have to choose
to connect to one of these dangling branches to get an parent. In
\cite{koskela-jenkins-johansen-spano-2020,mohle-1998}, the assumption
is that the probability of choosing one of the branches coming out
of particle 1 is 
\begin{equation}
\frac{\nu_{t}^{(1)}}{N}\label{eq:predicted-probability}
\end{equation}
 ($=2/4$ in our drawing). This is what makes possible to prove Equation
(\ref{eq:transition-formula}). Our goal is to show that this underlying
assumption does not hold under Hypothesis \ref{hyp:parental-distribution}.
\end{rem}

\section{Counterexample\label{sec:Counterexample}}

\subsection{Genealogical process in a simple case}

\begin{figure}[h]
\begin{centering}
\includegraphics{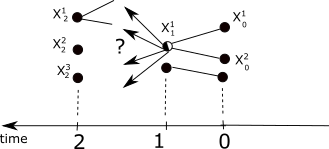}
\par\end{centering}
\caption{Simple case\label{fig:Simple-case}}
\end{figure}

We are interested in a very simple IPS. The state space is made of
two points: $a$ and $b$. The potential $g$ is such that $g(a)=p_{a}>0$,
$g(b)=p_{b}>0$ with 
\[
p_{b}<p_{a}\,.
\]
The terminal time is $T=2$ (the starting time is $0$). The law $\mu$
is such that $\mu(\{a\})=\alpha\in(0,1)$, $\mu(\{b\})=1-\alpha$.
The kernel $K$ is such that $K(\{a\},\{a\})=K(\{b\},\{b\})=1$. The
offspring count at time $T$ is such that
\[
\nu_{T}^{(1)}=2\,.
\]
Remember, this terminal time $T$ is the starting time of the IPS.
We suppose that $t=2$. We suppose that $\xi=\{\{1,2\},\dots\}$.
The situation is summarized in Figure \ref{fig:Simple-case}.
Each black dot represents a particle. For all $i$, $s$, $X_{s}^{(i)}$
is the position of the particle $i$ at time $s$ and we denote this
particle by \mybox[rounded corners]{mycol}{i,s}. The particle being
the parent of particles \mybox[rounded corners]{mycol}{1,0}, \mybox[rounded corners]{mycol}{2,0}
is the particle \mybox[rounded corners]{mycol}{1,1} in our drawing. 

In all the following, we will write $\boldsymbol{\nu}$ for $\boldsymbol{\nu}_{2}=(\nu_{2}^{(1)},\nu_{2}^{(2)},\dots,\nu_{2}^{(N)})$,
the offspring count at time $2$. This will lighten the notations
a little bit. 

\subsection{Idea of the counterexample}

We consider particle \mybox[rounded corners]{mycol}{1,1}. We know
that $\xi=\{\{1,2\},\dots\}$, that \mybox[rounded corners]{mycol}{1,1}
is the parent of \mybox[rounded corners]{mycol}{1,0} and \mybox[rounded corners]{mycol}{2,0},
that $\nu_{2}^{(1)}=2$. Given these informations (``\mybox[rounded corners]{mycol}{1,1}
has a lot of children'') and as $p_{a}>p_{b}$, $X_{1}^{(1)}$ is
more likely to be equal to $a$ than to $b$. But is can be equal
to $a$ only if the parent of \mybox[rounded corners]{mycol}{ 1,1}
is at position $a$. So the choice of a dangling branch (that will
induce the choice of a parent) might not be uniform among all those
dangling branches. 

This heuristics is rather poor and one could defend an opposite view
with arguments of the same type. This is why we resort to exact computation
in the next Section. 

\subsection{Main result}

We set
\[
f\,:\,x\in\R\mapsto\frac{x^{2}}{2}e^{-x}
\]
\[
q_{a}=\frac{p_{a}}{\alpha p_{a}+(1-\alpha)p_{b}}\,,\,q_{b}=\frac{p_{b}}{\alpha p_{a}+(1-\alpha)p_{b}}\,,
\]
\[
q_{a}'=\frac{p_{a}}{\alpha p_{a}q_{a}+(1-\alpha)p_{b}q_{b}}\,,\,q_{b}'=\frac{p_{b}}{\alpha p_{a}q_{a}+(1-\alpha)p_{b}q_{b}}\,.
\]

\begin{thm}
\label{th:main}We have (when $n$ is fixed)
\begin{multline}
\liminf_{N\rightarrow+\infty}\frac{N}{2}\times\p(a_{2}^{(1)}=1|\nu_{2}^{(1)}=2,\nu_{1}^{(1)}=2)\\
\geq\underset{=:R(\alpha,p_{a},p_{b})}{\underbrace{\frac{1}{\alpha q_{a}}\times\frac{\alpha f(q_{a})}{\alpha f(q_{a})+(1-\alpha)f(q_{a})}\times\frac{\alpha q_{a}f(q_{a}')}{\alpha q_{a}f(q_{a}')+(1-\alpha)q_{b}f(q_{b}')}}.}\label{eq:main-inequality}
\end{multline}
(We remind the reader that $a_{2}^{(1)}$ is the index of the parent
of particle \mybox[rounded corners]{mycol}{ j,1}.)
\end{thm}

\begin{proof}[Proof of Theorem \ref{th:main}]

We have
\begin{multline}
\p(a_{2}^{(1)}=1|\nu_{2}^{(1)}=2,\nu_{1}^{(1)}=2)\geq\p(X_{2}^{(1)}=a,X_{1}^{(1)}=a,a_{2}^{(1)}=1|\nu_{2}^{(1)}=2,\nu_{1}^{(1)}=2)\\
=\p(a_{2}^{(1)}=1|\nu_{2}^{(1)}=2,\nu_{1}^{(1)}=2,X_{2}^{(1)}=a,X_{1}^{(1)}=a)\times\p(X_{2}^{(1)}=a,X_{1}^{(1)}=a|\nu_{2}^{(1)}=2,\nu_{1}^{(1)}=2)\\
=\underset{=:(1)}{\underbrace{\p(a_{2}^{(1)}=1|\nu_{2}^{(1)}=2,\nu_{1}^{(1)}=2,X_{2}^{(1)}=a,X_{1}^{(1)}=a)}}\\
\times\underset{=:(2)}{\underbrace{\p(X_{2}^{(1)}=a|X_{1}^{(1)}=a,\nu_{2}^{(1)}=2,\nu_{1}^{(1)}=2)}}\times\underset{=:(3)}{\underbrace{\p(X_{1}^{(1)}=a|\nu_{2}^{(1)}=2,\nu_{1}^{(1)}=2)}}\label{eq:dec-aj2}
\end{multline}

We set, for $s=1,2$,
\[
N_{s}^{a}=\#\{k\,:\,X_{s}^{(k)}=a\}\,,\,N_{s}^{b}=\#\{k\,:\,X_{s}^{(k)}=b\}\,.
\]
As the only way of being equal to $a$ at time $1$ is to be the offspring
of a particle equal to $a$ at time $2$, we have 
\[
\sum_{k\,:\,X_{2}^{(k)}=a}\nu_{2}^{(k)}=N_{1}^{a}\,,
\]
and so
\[
(1)=\E\left(\left.\frac{\nu_{2}^{(1)}}{N_{1}^{a}}\right|\nu_{2}^{(1)}=2,\nu_{1}^{(1)}=2,X_{2}^{(1)}=a,X_{1}^{(1)}=a\right)\,.
\]
We set, for all $N\geq0$, 
\[
A_{N}=\frac{p_{a}N_{2}^{a}}{p_{a}N_{2}^{a}+p_{b}N_{2}^{b}}=\frac{p_{a}(N_{2}^{a}/N)}{p_{a}(N_{2}^{a}/N)+p_{b}(N_{2}^{b}/N)}.
\]
We have 
\[
A_{N}\convN\frac{p_{a}\alpha}{p_{a}\alpha+p_{b}(1-\alpha)}=:A_{\infty}\,.
\]
Furthermore, 
\[
x\in[0,1]\mapsto\frac{p_{a}x}{p_{a}x+p_{b}(1-x)}
\]
is $(p_{a}/p_{b})$-Lipschitz (because $p_{b}<p_{a}$), so, for all
$N$, 
\[
|A_{N}-A_{\infty}|\leq\frac{p_{a}}{p_{b}}\times\left|\frac{N_{2}^{a}}{N}-\alpha\right|\,,
\]
which implies (as $N_{2}^{a}$ is a sum of $N$ Bernoulli variables
of parameter $\alpha$)
\[
\E((A_{N}-A_{\infty})^{2})\leq\left(\frac{p_{a}}{p_{b}}\right)^{2}\times\frac{\alpha(1-\alpha)}{N}\,.
\]
We have 
\[
\frac{N_{1}^{a}}{N}=\frac{1}{N}\sum_{i=1}^{N}\1_{U_{i}\leq A_{N}}
\]
for some i.i.d. $U_{i}$'s, independent of the $A_{N}$'s (this is
what happens when we draw according to a categorical distribution).
By Lemma \ref{lem:If-we-have}, we then have
\begin{equation}
\frac{N_{1}^{a}}{N}\convN\frac{p_{a}\alpha}{p_{a}\alpha+p_{b}(1-\alpha)}=\alpha q_{a}\,,\label{eq:conv-Na1}
\end{equation}
which implies (as $N_{1}^{a}+N_{1}^{b}=N)$
\begin{equation}
\frac{N_{1}^{b}}{N}\convN\frac{p_{b}(1-\alpha)}{p_{a}\alpha+p_{b}(1-\alpha)}=(1-\alpha)q_{b}\,.\label{eq:conv-Nb1}
\end{equation}
So, by Fatou's Lemma, 
\begin{equation}
\liminf_{N\rightarrow+\infty}N\times(1)\geq2\times\frac{p_{a}\alpha+p_{b}(1-\alpha)}{p_{a}\alpha}=\frac{2}{\alpha q_{a}}\,.\label{eq:conv-(1)}
\end{equation}

We have 
\begin{eqnarray*}
(2) & = & \p(X_{2}^{(1)}=a|\nu_{2}^{(1)}=2,X_{1}^{(1)}=a,\nu_{1}^{(1)}=2)\\
 & = & \frac{\p(X_{2}^{(1)}=a,\nu_{2}^{(1)}=2,X_{1}^{(1)}=a)}{\p(\nu_{2}^{(1)}=2,X_{1}^{(1)}=a)}\,.
\end{eqnarray*}
Now, as a propagation of chaos consequence, we readily believe that
\mybox[rounded corners]{mycol}{1,2} and its descendants become independent
of \mybox[rounded corners]{mycol}{1,1} and its descendants when $N\rightarrow+\infty$
(see Section \ref{subsec:Proof-of-Equation} for a full proof). And
so
\begin{equation}
(2)=\frac{\p(X_{2}^{(1)}=a,\nu_{2}^{(1)}=2)\times\p(X_{1}^{(1)}=a,\nu_{1}^{(1)}=2)+O(1/\sqrt{N})}{\p(\nu_{2}^{(1)}=2)\times\p(X_{1}^{(1)}=a,\nu_{1}^{(1)}=2)+O(1/\sqrt{N})}\,.\label{eq:propagation-chaos-(2)}
\end{equation}
By the Law of Large Numbers,
\begin{equation}
\frac{N_{2}^{a}}{N}\convN\alpha\,,\,\frac{N_{2}^{b}}{N}\convN1-\alpha\,,\label{eq:LLN-01}
\end{equation}
 so we have (using the Dominated Convergence Theorem)
\begin{equation}
\p(X_{1}^{(1)}=a)=\E\left(\frac{p_{a}N_{2}^{a}}{p_{a}N_{2}^{a}+p_{b}N_{2}^{b}}\right)\underset{N\rightarrow+\infty}{\longrightarrow}\frac{\alpha p_{a}}{\alpha p_{a}+(1-\alpha)p_{b}}=\alpha q_{a}\,.\label{eq:lim-X-j-1-01}
\end{equation}
In the same way:
\begin{equation}
\p(X_{1}^{(1)}=b)\underset{N\rightarrow+\infty}{\longrightarrow}\frac{(1-\alpha)p_{b}}{\alpha p_{a}+(1-\alpha)p_{b}}=(1-\alpha)q_{b}\,.\label{eq:lim-X-j-1-02}
\end{equation}
We have 
\begin{multline*}
\p(\nu_{1}^{(1)}=2|X_{1}^{(1)}=a)=\E\left(\frac{N(N-1)}{2}\left(1-\frac{p_{a}}{N_{1}^{a}p_{a}+N_{1}^{b}p_{b}}\right)^{N-2}\left(\frac{p_{a}}{N_{1}^{a}p_{a}+N_{1}^{b}p_{b}}\right)^{2}\right)\\
=\E\left(\exp\left((N-2)\log\left(1-\frac{p_{a}}{N_{1}^{a}p_{a}+N_{1}^{b}p_{b}}\right)\right)\times\frac{N(N-1)}{2}\left(\frac{p_{a}}{N_{1}^{a}p_{a}+N_{1}^{b}p_{b}}\right)^{2}\right)\,.
\end{multline*}
Using Equations (\ref{eq:conv-Na1}), (\ref{eq:conv-Nb1}) and the
Dominated Convergence Theorem, we get
\begin{equation}
\p(\nu_{1}^{(1)}=2|X_{1}^{(1)}=a)\underset{N\rightarrow+\infty}{\longrightarrow}\frac{1}{2}\left(q_{a}'\right)^{2}e^{-q_{a}'}=f(q_{a}')\,.\label{eq:lim-nu-j-1-01}
\end{equation}
In the same way:
\begin{equation}
\p(\nu_{1}^{(1)}=2|X_{1}^{(1)}=b)\underset{N\rightarrow+\infty}{\longrightarrow}f(q_{b}')\,.\label{eq:lim-nu-j-1-02}
\end{equation}
So we can simplify the term (2) into
\begin{equation}
(2)=\frac{\p(X_{2}^{(1)}=a,\nu_{2}^{(1)}=2)+O(1/\sqrt{N})}{\p(\nu_{2}^{(1)}=2)+O(1/\sqrt{N})}\label{eq:propagation-chaos-01}
\end{equation}
Using Bayes Formula and the Formula of Total Probability, we get
\begin{eqnarray*}
\p(X_{2}^{(1)}=a|\nu_{2}^{(1)}=2) & = & \p(\nu_{2}^{(1)}=2|X_{2}^{(1)}=a)\frac{\p(X_{2}^{(1)}=a)}{\p(\nu_{2}^{(1)}=2)}\\
 & = & \frac{\p(\nu_{2}^{(1)}=2|X_{2}^{(1)}=a)\p(X_{2}^{(1)}=a)}{\p(\nu_{2}^{(1)}=2|X_{2}^{(1)}=a)\p(X_{2}^{(1)}=a)+\p(\nu_{2}^{(1)}=2|X_{2}^{(1)}=b)\p(X_{2}^{(1)}=b)}\\
 & = & \frac{\p(\nu_{2}^{(1)}=2|X_{2}^{(1)}=a)\times\alpha}{\p(\nu_{2}^{(1)}=2|X_{2}^{(1)}=a)\times\alpha+\p(\nu_{2}^{(1)}=2|X_{2}^{(1)}=b)\times(1-\alpha)}\,.
\end{eqnarray*}
We have 
\begin{multline*}
\p(\nu_{2}^{(1)}=2|X_{2}^{(1)}=a)=\E\left(\frac{N(N-1)}{2}\left(1-\frac{p_{a}}{N_{2}^{a}p_{a}+N_{2}^{b}p_{b}}\right)^{N-2}\left(\frac{p_{a}}{N_{2}^{a}p_{a}+N_{2}^{b}p_{b}}\right)^{2}\right)\\
=\E\left(\exp\left((N-2)\log\left(1-\frac{p_{a}}{N_{2}^{a}p_{a}+N_{2}^{b}p_{b}}\right)\right)\times\frac{N(N-1)}{2}\left(\frac{p_{a}}{N_{2}^{a}p_{a}+N_{2}^{b}p_{b}}\right)^{2}\right)\,.
\end{multline*}
Using Equation (\ref{eq:LLN-01}) and the Dominated Convergence Theorem,
we get
\begin{equation}
\p(\nu_{2}^{(1)}=2|X_{2}^{(1)}=a)\times\alpha\underset{N\rightarrow+\infty}{\longrightarrow}\frac{\alpha}{2}q_{a}^{2}e^{-q_{a}}=\alpha f(q_{a})\,.\label{eq:lim-nu-1-2-01}
\end{equation}
In the same way, we get
\begin{equation}
\p(\nu_{2}^{(1)}=2|X_{2}^{(1)}=b)\times(1-\alpha)\underset{N\rightarrow+\infty}{\longrightarrow}\frac{(1-\alpha)}{2}q_{b}^{2}e^{-q_{b}}=(1-\alpha)f(q_{b})\,.\label{eq:lim-nu-1-2-02}
\end{equation}
So we get 
\begin{equation}
(2)\underset{N\rightarrow+\infty}{\longrightarrow}\frac{\alpha f(q_{a})}{\alpha f(q_{a})+(1-\alpha)f(q_{b})}\,.\label{eq:conv-(2)}
\end{equation}

We have 
\[
(3)=\frac{\p(X_{1}^{(1)}=a,\nu_{2}^{(1)}=2,\nu_{2}^{(1)}=2)}{\p(\nu_{2}^{(1)}=2,\nu_{2}^{(1)}=2)}\,.
\]
Again, as a propagation of chaos consequence, we readily believe that
\mybox[rounded corners]{mycol}{1,2} and its descendants become independent
of \mybox[rounded corners]{mycol}{1,1} and its descendants when $N\rightarrow+\infty$.
And so
\begin{eqnarray*}
(3) & = & \frac{\p(X_{1}^{(1)}=a,\nu_{1}^{(1)}=2)\p(\nu_{2}^{(1)}=2)+O(1/\sqrt{N})}{\p(\nu_{2}^{(1)}=2,)\p(\nu_{2}^{(1)}=2)+O(1/\sqrt{N})}\\
 &  & \text{(by Equations (\ref{eq:LLN-01}), (\ref{eq:lim-nu-1-2-01}), (\ref{eq:lim-nu-1-2-02}), we can symplify)}\\
 & = & \frac{\p(X_{1}^{(1)}=a,\nu_{1}^{(1)}=2)+O(1/\sqrt{N})}{\p(\nu_{2}^{(1)}=2)+O(1/\sqrt{N})}
\end{eqnarray*}
(the proof the above Equation is very similar to the proof of Equation
(\ref{eq:propagation-chaos-01}) so we omit it). Using Bayes Formula
and the Formula of Total Probability, 
\begin{eqnarray*}
\p(X_{1}^{(1)}=a|\nu_{1}^{(1)}=2) & = & \p(\nu_{1}^{(1)}=2|X_{1}^{(1)}=a)\times\frac{\p(X_{1}^{(1)}=a)}{\p(\nu_{1}^{(1)}=2)}\\
 & = & \frac{\p(\nu_{1}^{(1)}=2|X_{1}^{(1)}=a)\p(X_{1}^{(1)}=a)}{\p(\nu_{1}^{(1)}=2|X_{1}^{(1)}=a)\p(X_{1}^{(1)}=a)+\p(\nu_{1}^{(1)}=2|X_{1}^{(1)}=b)\p(X_{1}^{(1)}=b)}\,.
\end{eqnarray*}
By Equations (\ref{eq:lim-X-j-1-01}), (\ref{eq:lim-X-j-1-02}), (\ref{eq:lim-nu-j-1-01}),
(\ref{eq:lim-nu-j-1-02}), we get

\begin{equation}
(3)\underset{N\rightarrow+\infty}{\longrightarrow}\frac{\alpha q_{a}f(q_{a}')}{\alpha q_{a}f(q_{a}')+(1-\alpha)q_{b}f(q_{b'})}\,.\label{eq:conv-(3)}
\end{equation}

From Equations (\ref{eq:dec-aj2}), (\ref{eq:conv-(1)}), (\ref{eq:conv-(2)}),
(\ref{eq:conv-(3)}), we get 
\begin{multline*}
\liminf_{N\rightarrow+\infty}N\p(a_{2}^{(1)}=1|\nu_{2}^{(1)}=2,\nu_{1}^{(1)}=2)\geq2\times\frac{1}{\alpha q_{a}}\times\frac{\alpha f(q_{a})}{\alpha f(q_{a})+(1-\alpha)f(q_{b})}\\
\times\frac{\alpha q_{a}f(q_{a}')}{\alpha q_{a}f(q_{a}')+(1-\alpha)q_{b}f(q_{b'})}
\end{multline*}

\end{proof}
\begin{cor}
\label{cor:For-well-chosen}For well chosen values of $\alpha$, $p_{a}$,
$p_{b}$, the right-hand $R(\alpha,p_{a},p_{b})$ side of inequality
\ref{eq:main-inequality} in Theorem \ref{th:main} is strictly bigger
than $1.$This contradicts the claim of Equation (\ref{eq:predicted-probability}).
\end{cor}

\begin{proof}
The trick is to choose the values of $\alpha$, $p_{a}$, $p_{b}$
and to compute $R(\alpha,p_{a},p_{b})$. We present a \texttt{python}
code for computing $R$ below. We begin by the definition of $q_{a}$
and $q_{b}$. 

\begin{lstlisting}
import numpy as np 
import matplotlib.pyplot as plt
def qa(al,pa,pb):     
	return(pa/(al*pa+(1-al)*pb)) 
def qb(al,pa,pb):     
	return(pb/(al*pa+(1-al)*pb))
\end{lstlisting}
Here we define $q_{a}'$, $q_{b}'$.
\begin{lstlisting}
def qap(al,pa,pb):     
	return(pa/(al*qa(al,pa,pb)*pa+(1-al)*qb(al,pa,pb)*pb)) 
def qbp(al,pa,pb):     
	return(pb/(al*qa(al,pa,pb)*pa+(1-al)*qb(al,pa,pb)*pb))
\end{lstlisting}
Here, we define $f$ and the first term ($1/(\alpha q_{a})$) and
the second term ($\alpha f(q_{a})/(\alpha f(q_{a})+(1-\alpha)f(q_{b}))$)
appearing in $R(\alpha,p_{a},p_{b})$.
\begin{lstlisting}
def f(t):    
	return(0.5*t**2*np.exp(-t)) 
def t1(al,pa,pb):     
	return(1/(al*qa(al,pa,pb))) 
def t2(al,pa,pb):    
	return((al*f(qa(al,pa,pb)))/
		(al*f(qa(al,pa,pb))+(1-al)*f(qb(al,pa,pb))))
\end{lstlisting}
Here, we define the third term ($\alpha q_{a}f(q_{a}')/(\alpha q_{a}f(q_{a}')+(1-\alpha)q_{b}f(q_{b}'))$)
appearing in $R(\alpha,p_{a},p_{b})$.
\begin{lstlisting}
def t3(al,pa,pb):     
	return(al*qa(al,pa,pb)*f(qap(al,pa,pb))/
		(al*qa(al,pa,pb)*f(qap(al,pa,pb))
		+(1-al)*qb(al,pa,pb))*f(qbp(al,pa,pb))))
\end{lstlisting}
Now, $R$ is defined by the following.
\begin{lstlisting}
def R(al,pa,pb):    
	return(t1(al,pa,pb)*t2(al,pa,pb)*t3(al,pa,pb))
\end{lstlisting}
We draw $R(0.5,1,p_{b})$ for $p_{b}$ in $[0,0.2]$.
\begin{lstlisting}
def trace(x):    
	return(R(0.5,1,x))
x=np.linspace(0,0.2,500) 
plt.plot(x,trace(x))
plt.xlabel('x') 
plt.ylabel('R(0.5,1,x)')
plt.show()
\end{lstlisting}
We get Figure \ref{fig:The-term-.}.
\begin{figure}[h]
\begin{centering}
\includegraphics[scale=0.75]{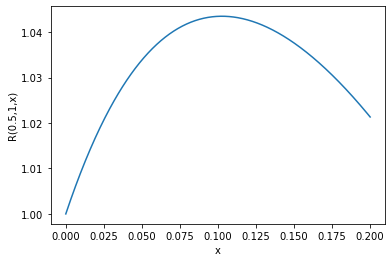}
\par\end{centering}
\caption{\label{fig:The-term-.}The term $R$. }
\end{figure}
The line above $1.$ Furthermore, if we ask for $R(0.5,1,0.075)$,
we get the answer $1.04104942\dots$ which is strictly bigger than
$1$. 

We have
\begin{multline*}
\E(\p(a_{2}^{(1)}=1|\boldsymbol{\nu},\xi)|\nu_{2}^{(1)}=2,\nu_{1}^{(1)}=2)=\E(\E(\1_{a_{2}^{(1)}=1}|\boldsymbol{\nu},\xi)|\nu_{2}^{(1)}=2,\nu_{1}^{(1)}=2))\\
=\E(\1_{a_{2}^{(1)}=1}|\nu_{2}^{(1)}=2,\nu_{1}^{(1)}=2)=\p(a_{2}^{(1)}=1|\nu_{2}^{(1)}=2,\nu_{1}^{(1)}=2)
\end{multline*}
So, if we take $\alpha=0.5$, $p_{a}=1$, $p_{b}=0.075$, there exist
$N$ and $\boldsymbol{\nu}$, $\xi$ (with $\nu_{2}^{(1)}=2$, $\nu_{1}^{(1)}=2$)
such that 
\[
\p(a_{2}^{(1)}=1|\boldsymbol{\nu},\xi)>\frac{\nu_{2}^{(1)}}{N}.
\]
This contradicts Equation (\ref{eq:predicted-probability}).
\end{proof}

\section{Appendix}

\subsection{Technical probability Lemma }
\begin{lem}
\label{lem:If-we-have}If we have a sequence of $L^{2}$ random variables
$(A_{N})_{N\geq0}$ in $[0,1]$ converging a.s. towards a constant
$A_{\infty}$ as $N\rightarrow+\infty$ and
\[
\E((A_{N}-A_{\infty})^{2})\leq\frac{C}{N}
\]
(for some constant $C$) and if we have an i.i.d. sequence $(U_{i})_{i\geq0}$
of variables of law $\mathcal{U}([0,1])$ (independent of the $A_{n}$'s)
then
\[
\frac{1}{N}\sum_{i=1}^{N}\1_{U_{i}\leq A_{N}}\overset{\text{a.s.}}{\underset{N\rightarrow+\infty}{\longrightarrow}}A_{\infty}\,.
\]
\end{lem}

\begin{proof}
We have (law of large numbers)
\[
\frac{1}{N}\sum_{i=1}^{N}\1_{U_{i}\leq A_{\infty}}\overset{\text{a.s.}}{\underset{N\rightarrow+\infty}{\longrightarrow}}A_{\infty}\,.
\]
So we look at the difference
\begin{multline*}
\left|\frac{1}{N}\sum_{i=1}^{N}\1_{U_{i}\leq A_{N}}-\frac{1}{N}\sum_{i=1}^{N}\1_{U_{i}\leq A_{\infty}}\right|=\left|\frac{1}{N}\sum_{i=1}^{N}u(U_{i},A_{N},A_{\infty})\right|\\
\leq\left|\frac{1}{N}\sum_{i=1}^{N}u(U_{i},A_{N},A_{\infty})-(A_{N}-A_{\infty})\right|+\frac{1}{N}\sum_{i=1}^{N}|A_{N}-A_{\infty}|\,,
\end{multline*}
where
\[
u\,:\,(v,x,y)\in[0,1]\times\R\times\R\mapsto u(v,x,y)=\begin{cases}
-1 & \text{ if }x<y\text{ and }v\in(x,y]\text{,}\\
1 & \text{ if }y<x\text{ and }v\in(y,x]\text{,}\\
0 & \text{ otherwise.}
\end{cases}
\]
We have 
\[
\frac{1}{N}\sum_{i=1}^{N}|A_{N}-A_{\infty}|=|A_{N}-A_{\infty}|\convN0\,.
\]
And we have (as $\E(u(U_{i},A_{N},A_{\infty})|A_{N},A_{\infty})=A_{N}-A_{\infty}$
for all $i\geq0$, $N\geq0$)
\begin{multline*}
\E\left(\left(\frac{1}{N}\sum_{i=1}^{N}u(U_{i},A_{N},A_{\infty})-(A_{N}-A_{\infty})\right)^{2}\right)=\E\left(\frac{1}{N^{2}}\sum_{i=1}^{N}(u(U_{i},A_{N},A_{\infty})-(A_{N}-A_{\infty}))^{2}\right)\\
=\E\left(\frac{1}{N^{2}}\sum_{i=1}^{N}\E((u(U_{i},A_{N},A_{\infty})-(A_{N}-A_{\infty}))^{2}|A_{N},A_{\infty})\right)\\
=\E\left(\frac{1}{N^{2}}\sum_{i=1}^{N}\E(u(U_{i},A_{N},A_{\infty})^{2}|A_{N},A_{\infty})-(A_{N}-A_{\infty})^{2}\right)\\
=\E\left(\frac{1}{N^{2}}\sum_{i=1}^{N}|A_{N}-A_{\infty}|-(A_{N}-A_{\infty})^{2}\right)\\
\leq\E\left(\frac{1}{N^{2}}\sum_{i=1}^{N}|A_{N}-A_{\infty}|\right)\\
\text{(Jensen's inequality)}\leq\frac{N}{N^{2}}\E((A_{N}-A_{\infty})^{2})^{1/2}\leq\frac{\sqrt{C}}{N^{3/2}}\,.
\end{multline*}
By Borell-Cantelli's Lemma:
\[
\frac{1}{N}\sum_{i=1}^{N}u(U_{i},A_{N},A_{\infty})-(A_{N}-A_{\infty})\convN0\,,
\]
and this finishes our proof.
\end{proof}

\subsection{\label{subsec:Proof-of-Equation}Proof of Equation (\ref{eq:propagation-chaos-(2)})}

To prove Equation (\ref{eq:propagation-chaos-(2)}), it is sufficient
to prove the following Lemma.
\begin{lem}
\label{lem:propagation-of-chaos}We have
\begin{equation}
\p(X_{2}^{(1)}=a,\nu_{2}^{(1)}=2,X_{1}^{(1)}=a,\nu_{1}^{(1)}=2)=\p(X_{2}^{(1)}=a,\nu_{2}^{(1)}=2)\p(X_{1}^{(1)}=a,\nu_{1}^{(1)}=2)+O(1/\sqrt{N})\,,\label{eq:p-c-01}
\end{equation}
\end{lem}

\begin{equation}
\p(\nu_{2}^{(1)}=1,X_{1}^{(1)}=a,\nu_{1}^{(1)}=2)=\p(\nu_{2}^{(1)}=1)\p(X_{1}^{(1)}=a,\nu_{1}^{(1)}=2)+O(1/\sqrt{N})\,.\label{eq:p-c-02}
\end{equation}
The proof of the above Lemma will be done through coupling.

\subsubsection{Description of our coupling}

When sampling from the categorical distribution with weights $(w_{2}^{(1)},\dots,w_{2}^{(N)}$)
(these are the weights of particles at time $2$), we suppose we use
i.i.d. variables $(U_{2}^{(1)},\dots,U_{2}^{(N)})$ of law $\mathcal{U}([0,1])$
in the following way
\[
a_{2}^{(i)}=k\text{ if and only if }w_{2}^{(1)}+\dots+w_{2}^{(k-1)}\leq U_{2}^{(i)}<w_{2}^{(1)}+\dots+w_{2}^{(k)}\text{ (}k\in[N]\text{).}
\]
We remind the reader that, for all $i$ in $[N]$,
\[
w_{2}^{(i)}=\frac{g(X_{2}^{(i)})}{g(X_{2}^{(1)})+\dots+g(X_{2}^{(N)})}\,.
\]
We introduce 
\[
\widehat{w}_{2}^{(1)}=\frac{g(X_{2}^{(1)})}{N(\alpha p_{a}+(1-\alpha)p_{b})}\,,
\]
\[
\widehat{\nu}_{2}^{(1)}=\#\{i\,:\,U_{2}^{(i)}\leq\widehat{w}_{2}^{(1)}\}\,.
\]
The random number $\widehat{\nu}_{2}^{(1)}$ plays the role of a ``number
of descendants'' of \mybox[rounded corners]{mycol}{1,2} which is
independent of $X_{2}^{(2)},\dots,X_{2}^{(N)}$. 

We introduce additional i.i.d. variables $(U_{2}^{(N+1)},U_{2}^{(N+2)},\dots$)
of law $\mathcal{U}([0,1])$. %
\begin{comment}
We also introduce i.i.d. variables $(X_{2}^{(N+1)},X_{2}^{(N+2)},\dots)$
of law $\mu$.
\end{comment}
{} We then set, for all $i$ in $\{N+1,N+2,\dots\}$, $k$ in $[N]$,
\[
a_{_{2}}^{(i)}=k\text{ if and only if \ensuremath{w_{2}^{(1)}}+\ensuremath{\dots}+\ensuremath{w_{2}^{(k-1)}\leq U_{2}^{(i)}}}<w_{2}^{(1)}+\dots+w_{2}^{(k)}\,.
\]
For all $i$ in $\{1,2,\dots,N,N+1,\dots\}$, we set 
\[
X_{1}^{(i)}=X_{2}^{a_{2}^{(i)}}
\]
(with our overly simple $K$, this amounts to sampling $X_{1}^{(i)}\sim K(X_{2}^{a_{2}^{(i)}},.)$).%
\begin{comment}
. We observe that $(X_{1}^{(1)},\dots,X_{1}^{(N)})$ are (in law)
identical to the $(X_{1}^{(1)},\dots,X_{1}^{(N)})$ described in Algorithm
\ref{alg:Simulation-of-an}, in the case of our simplified model.
\end{comment}
{} %
\begin{comment}
We have added a random number of particles of values 
\[
(X_{1}^{(N+1)},\dots,X_{1}^{(N+\nu_{2}^{(1)})})
\]
 in $E$. For each of these particle, its parent is not \mybox[rounded corners]{mycol}{1,2}.
We observe, that amongst \mybox[rounded corners]{mycol}{1,1} , ...
, \mybox[rounded corners]{mycol}{\nu^{(1)}_2,1} , exactly $N$ particles
have \mybox[rounded corners]{mycol}{1,2} for parent. We define a
random set: 
\[
\mathcal{I}=\{i\in[N+\nu_{2}^{(1)}]\,:\,a_{2}^{(i)}\neq1\}.
\]
 We write 
\end{comment}

We define $\mathcal{I}=\{i_{1},i_{2},\dots,i_{N}\}$ with $i_{1}<i_{2}<\dots<i_{N}$
in a recursive way:
\begin{itemize}
\item $i_{1}=\inf\{i\,:\,a_{2}^{(i)}\neq1\}$
\item for $k\in\{2,\dots,N\}$, $i_{k}=\inf\{i>i_{k-1}\,:\,a_{2}^{(i_{k})}\neq1\}$.
\end{itemize}
The $N$-uple 
\begin{equation}
(X_{1}^{(i_{1})},\dots,X_{1}^{(i_{N})})=(X_{2}^{(a_{2}^{(i_{1})})},\dots,X_{2}^{(a_{2}^{(i_{N})})})\text{ is independent of }(X_{2}^{(1)},\widehat{\nu}_{2}^{(1)})\,.\label{eq:I-independant}
\end{equation}
Indeed, the sequence $i_{1},i_{2},\dots$ is built with an accep-reject
scheme: each $U_{2}^{(i_{k})}$ is of law $\mathcal{U}([w_{2}^{(1)},1])$
and for all $j$ in $\{2,\dots N\}$, $a_{2}^{(i_{k})}=j$ with probability
\[
\frac{g(X_{2}^{(j)})}{\sum_{r=2}^{N}g(X_{2}^{(r)})}\,.
\]
 We set 
\[
I=i_{1}\,,
\]
and 
\[
\widetilde{X}_{1}^{(1)}=X_{1}^{(I)}\,.
\]
\begin{comment}
The variable $\widetilde{X}_{1}^{(1)}$ is independent of $(X_{2}^{(1)},\nu_{2}^{(1)})$. 
\end{comment}

Suppose we use i.i.d.. variables $(U_{1}^{(1)},\dots,U_{1}^{(N)})$
of law $\mathcal{U}([0,1])$ to compute the parental indices (for
all $i$ in $[N])$:
\[
a_{1}^{(i)}=k\text{ if and only if }w_{1}^{(1)}+\dots+w_{1}^{(k-1)}\leq U_{1}^{(i)}<w_{1}^{(1)}+\dots+w_{1}^{(k)}\text{ (}k\in[N]\text{)},
\]
where the $(w_{1}^{(1)},\dots,w_{1}^{(N)})$ are the weighs of the
particle at time $1$. For all $i$ in $[N]$, 
\[
w_{1}^{(i)}=\frac{g(X_{1}^{(i)})}{\sum_{r=1}^{N}g(X_{1}^{(r)})}\,.
\]
We introduce alternative weights
\[
\widetilde{w}_{1}^{(i)}=\frac{g(X_{1}^{(i)})}{\sum_{k\in\mathcal{I}}g(X_{1}^{(k)})}\text{ (}i\in\mathcal{I}\text{).}
\]
From these, we compute alternative parental indices (for all $i$,
$k$ in $[N]$):
\[
\widetilde{a}_{1}^{(i)}=k\text{ if and only if }\widetilde{w}_{1}^{(i_{1})}+\dots+\widetilde{w}_{1}^{(i_{k})}\leq U_{1}^{(i)}\leq\widetilde{w}_{1}^{(i_{1})}+\dots+\widetilde{w}_{1}^{(i_{k})}\,.
\]
And we set 
\[
\widetilde{\nu}_{1}^{(1)}=\#\{i\in[N]\,:\,\widetilde{a}_{1}^{(i)}=I\}\,.
\]

\subsubsection{Technical Lemmas}
\begin{lem}
\label{lem:independent}The variable $(\widetilde{X}_{1}^{(1)},\widetilde{\nu}_{1}^{(1)})$
is independent of $(X_{2}^{(1)},\widehat{\nu}_{2}^{(1)})$.
\end{lem}

\begin{proof}
We want to show, that for all $u\in\{a,b\}$ and $k\in\N$, $\p((\widetilde{X}_{1}^{(1)},\widetilde{\nu}_{1}^{(1)})=(u,k)|X_{2}^{(1)},\widehat{\nu}_{2}^{(1)})$
is a constant. We set 
\[
\beta(u)=\begin{cases}
\alpha & \text{ if }u=a,\\
1-\alpha & \text{ if }u=b\,.
\end{cases}
\]
 We have, for all $u$ in $\{a,b\}$,

\begin{multline*}
\p(\widetilde{X}_{1}^{(1)}=u|X_{2}^{(1)},\widehat{\nu}_{2}^{(1)})=\sum_{j=2}^{N}\p(X_{2}^{(j)}=u,a_{2}^{I}=j|X_{2}^{(1)},\widehat{\nu}_{2}^{(1)})\\
=\sum_{j=2}^{N}\p(a_{2}^{I}=j|X_{2}^{(1)},\widehat{\nu}_{2}^{(1)},X_{2}^{(j)}=u)\p(X_{2}^{(j)}=u|X_{2}^{(1)},\widehat{\nu}_{2}^{(1)})\\
=\sum_{j=2}^{N}\E\left(\E\left.\left(\left.\1_{a_{2}^{I}=j}\right|X_{2}^{(1)},\widehat{\nu}_{2}^{(1)},(X_{2}^{(i)})_{2\leq i\leq N}\right)\right|X_{2}^{(1)},\widehat{\nu}_{2}^{(1)},X_{2}^{(j)}=u\right)\times\p(X_{2}^{(j)}=u|X_{2}^{(1)},\widehat{\nu}_{2}^{(1)})\\
=\sum_{j=2}^{N}\E\left(\left.\frac{g(u)}{\sum_{i=2}^{N}g(X_{2}^{(i)})}\right|X_{2}^{(1)},\widehat{\nu}_{2}^{(1)},X_{2}^{(j)}=u\right)\times\beta(u)\\
=\sum_{j=2}^{N}\E\left(\frac{g(u)}{g(u)+\sum_{i\in\{2,\dots,N\}\backslash\{j\}}g(X_{1}^{(i)})}\right)\times\beta(u)\,.
\end{multline*}

So $\p(\widetilde{X}_{1}^{(1)}=u|X_{2}^{(1)},\widehat{\nu}_{2}^{(1)})$
is constant. 

Next, we have (because of Equation (\ref{eq:I-independant})), for
all $k$ in $\{0,1,\dots,N\}$, 
\begin{eqnarray*}
\p(\widetilde{\nu}_{1}^{(1)}=k|\mathcal{I},(X_{1}^{(i)})_{i\in\mathcal{I}},X_{2}^{(1)},\widehat{\nu}_{2}^{(1)}) & = & \binom{N}{k}(\widetilde{w}_{1}^{(I)})^{k}(1-\widetilde{w}_{1}^{(I)})^{N-k}\\
 & = & \binom{N}{k}\left(\frac{g(\widetilde{X}_{1}^{(1)})}{\sum_{i\in\mathcal{I}}g(X_{1}^{(i)})}\right)^{k}\left(1-\frac{g(\widetilde{X}_{1}^{(1)})}{\sum_{i\in\mathcal{I}}g(X_{1}^{(i)})}\right)^{N-k}\,.
\end{eqnarray*}
 So 
\begin{multline*}
\p(\widetilde{\nu}_{1}^{(1)}=k|\widetilde{X}_{1}^{(1)}=u,X_{2}^{(1)},\widehat{\nu}_{2}^{(1)})\\
=\E\left(\binom{N}{k}\left(\frac{g(u)}{g(u)+\sum_{i\in\mathcal{I}\backslash\{i_{1}\}}g(X_{1}^{(i)})}\right)^{k}\left(1-\frac{g(u)}{g(u)+\sum_{i\in\mathcal{I}\backslash\{i_{1}\}}g(X_{1}^{(i)})}\right)^{N-k}\right)
\end{multline*}
 is constant.

So we have that, for all $k$,
\[
\p(\widetilde{X}_{1}^{(1)}=u,\widetilde{\nu}_{1}^{(1)}=k|X_{2}^{(1)},\widehat{\nu}_{2}^{(1)})=\p(\widetilde{\nu}_{1}^{(1)}=k|X_{2}^{(1)},\widehat{\nu}_{2}^{(1)},\widetilde{X}_{1}^{(1)}=u)\p(\widetilde{X}_{1}^{(1)}=u|X_{2}^{(1)},\widehat{\nu}_{2}^{(1)})
\]
is a constant. And this finishes the proof. 
\end{proof}
\begin{lem}
\label{lem:O(1/N)}$\p((\widetilde{X}_{1}^{(1)},\widetilde{\nu}_{1}^{(1)})\neq(X_{1}^{(1)},\nu_{1}^{(1)}))=O(1/N)$
\end{lem}

\begin{proof}
We have
\begin{multline*}
\p((\widetilde{X}_{1}^{(1)},\widetilde{\nu}_{1}^{(1)})=(X_{1}^{(1)},\nu_{1}^{(1)}))\geq\p(I=1,\widetilde{\nu}_{1}^{(1)}=\nu_{1}^{(1)})=\p(a_{2}^{(1)}\neq1,\widetilde{\nu}_{1}^{(1)}=\nu_{1}^{(1)})\\
=\p(a_{2}^{(1)}\neq1)\p(\nexists i\in[N]\,:\,U_{1}^{(i)}\in[w_{1}^{(1)},\widetilde{w}_{1}^{(1)}]|a_{2}^{(1)}\neq1)\\
=\E(1-w_{2}^{(1)})\times\E\left(\left.(1-|w_{1}^{(1)}-\widetilde{w}_{1}^{(1)}|)^{N}\right|a_{2}^{(1)}\neq1\right)\text{}\\
\text{(for }N\text{ such that }\frac{p_{a}}{Np_{b}}<1\text{) }\\
\geq\left(1-\frac{p_{a}}{Np_{b}}\right)\times\E\left(\left.1-N|w_{1}^{(1)}-\widetilde{w}_{1}^{(1)}|\right|a_{2}^{(1)}\neq1\right)\\
\text{(because }a_{2}^{(1)}\neq1\Rightarrow i_{1}=1\text{) }\\
=\left(1-\frac{p_{a}}{Np_{b}}\right)\E\left(\left.1-N\left|\frac{g(X_{1}^{(1)})}{\sum_{i\in[N]}g(X_{1}^{(i)})}-\frac{g(X_{1}^{(1)})}{\sum_{i\in\mathcal{I}}g(X_{1}^{(i)})}\right|\right|a_{2}^{(1)}\neq1\right)\\
=\left(1-\frac{p_{a}}{Np_{b}}\right)\E\left(\left.1-N\left|\frac{g(X_{1}^{(1)})(\sum_{i\in[N]}g(X_{1}^{(i)})-\sum_{i\in\mathcal{I}}g(X_{1}^{(i)}))}{\sum_{i\in[N]}g(X_{1}^{(i)})\times\sum_{i\in\mathcal{I}}g(X_{1}^{(i)})}\right|\right|a_{2}^{(1)}\neq1\right)\\
\text{(as }\#([N]\backslash\mathcal{I})=\#(\mathcal{I}\backslash[N])=\nu_{1}^{(2)}\text{)}\geq\left(1-\frac{p_{a}}{Np_{b}}\right)\E\left(\left.1-\frac{\nu_{2}^{(1)}p_{a}(p_{a}-p_{b})}{Np_{b}^{2}}\right|a_{2}^{(1)}\neq1\right)\\
=\left(1-\frac{p_{a}}{Np_{b}}\right)\left(1-\frac{\E(\nu_{2}^{(1)}|a_{2}^{(1)}\neq1)p_{a}^{2}}{Np_{b}^{2}}\right)\\
\text{(as }\nu_{2}^{(1)}=\sum_{i=1}^{N}\1_{\{1\}}(a_{2}^{(i)})\text{)}\geq\left(1-\frac{p_{a}}{Np_{b}}\right)\left(1-\frac{\E(\nu_{2}^{(1)})p_{a}^{2}}{Np_{b}^{2}}\right)\\
\geq\left(1-\frac{p_{a}}{Np_{b}}\right)\left(1-\frac{\left(\frac{p_{a}}{p_{b}}\right)\times p_{a}^{2}}{Np_{b}^{2}}\right)
\end{multline*}
\end{proof}
\begin{lem}
\label{lem:O(1/sqrt(N))}$\p((X_{2}^{(1)},\widehat{\nu}_{2}^{(1)})\neq(X_{2}^{(1)},\nu_{2}^{(1)}))=O(1/\sqrt{N})$ 
\end{lem}

\begin{proof}
We have
\begin{multline*}
\p((X_{2}^{(1)},\widehat{\nu}_{2}^{(1)})\neq(X_{2}^{(1)},\nu_{2}^{(1)})=\p(\exists i\in[N]\,:\,U_{2}^{(i)}\in[w_{2}^{(1)},\widehat{w}_{2}^{(1)}])\\
=1-\p\left(\forall i\in[N],\,U_{2}^{(i)}\notin\left[\frac{g(X_{2}^{(1)})}{\sum_{i=1}^{N}g(X_{2}^{(i)})},\frac{g(X_{2}^{(1)})}{N(\alpha p_{a}+(1-\alpha)p_{b})}\right]\right)\\
=1-\E\left(\left(1-\left|\frac{g(X_{2}^{(1)})}{N(\alpha p_{a}+(1-\alpha)p_{b})}-\frac{g(X_{2}^{(1)})}{\sum_{i=1}^{N}g(X_{2}^{(i)})}\right|\right)^{N}\right)\\
\leq N\E\left(\left|\frac{g(X_{2}^{(1)})}{N(\alpha p_{a}+(1-\alpha)p_{b})}-\frac{g(X_{2}^{(1)})}{\sum_{i=1}^{N}g(X_{2}^{(i)})}\right|\right)\\
=N\E\left(\left|g(X_{2}^{(1)})\times\frac{\sum_{i=1}^{N}g(X_{2}^{(i)})-N(\alpha p_{a}+(1-\alpha)p_{b})}{N(\alpha p_{a}+(1-\alpha)p_{b})\times\sum_{i=1}^{N}g(X_{2}^{(i)})}\right|\right)\\
\text{(Jensen's inequality)}\leq\frac{Np_{a}}{N^{2}p_{b}^{2}}\E\left(\left(\sum_{i=1}^{N}(g(X_{2}^{(i)})-(\alpha p_{a}+(1-\alpha)p_{b})\right)^{2}\right)^{1/2}\\
\text{(with }\sigma\geq0\,,\,\sigma^{2}=\alpha p_{a}^{2}+(1-\alpha)p_{b}^{2}-(\alpha p_{a}+(1-\alpha)p_{b})^{2}\text{) }=\frac{p_{a}}{Np_{b}^{2}}\sqrt{N}\sigma
\end{multline*}
\end{proof}

\subsection{Proof of Lemma \ref{lem:propagation-of-chaos}}
\begin{proof}
We compute 
\begin{multline*}
\p(X_{2}^{(1)}=a,\nu_{2}^{(1)}=2,X_{1}^{(1)}=a,\nu_{1}^{(1)}=2)\\
\text{(by Lemma \ref{lem:O(1/N)} and Lemma \ref{lem:O(1/sqrt(N))})}\\
=\p(X_{2}^{(1)}=a,\widehat{\nu}_{2}^{(1)}=2,\widetilde{X}_{1}^{(1)}=a,\widetilde{\nu}_{1}^{(1)}=2)+O(1/\sqrt{N})\\
\text{(by Lemma \ref{lem:independent})}\\
=\p(X_{2}^{(1)}=a,\widehat{\nu}_{2}^{(1)}=2)\times\p(\widetilde{X}_{1}^{(1)}=a,\widetilde{\nu}_{1}^{(1)}=2)+O(1/\sqrt{N})\\
\text{(by Lemma \ref{lem:O(1/N)} and Lemma \ref{lem:O(1/sqrt(N))})}\\
=\p(X_{2}^{(1)}=a,\nu_{2}^{(1)}=2)\times\p(X_{1}^{(1)}=a,\nu_{1}^{(1)}=2)+O(1/\sqrt{N})
\end{multline*}
This proves Equation (\ref{eq:p-c-01}). The proof of Equation (\ref{eq:p-c-02})
is very similar. 
\end{proof}
\bibliographystyle{amsalpha}
\bibliography{coalescence-formula}

\end{document}